\documentclass[a4paper,leqno]{article}
\usepackage{latexsym,amssymb,amsmath}
\usepackage{color,amsfonts}
\usepackage[english]{babel}

\newtheorem{theorem}{Theorem}[section]
\newtheorem{lemma}[theorem]{Lemma}
\newtheorem{proposition}[theorem]{Proposition}

\newtheorem{definition}[theorem]{Definition}

\newenvironment{proof}[1][Proof]{\textbf{#1.} }{\ \mathbb Rule{0.5em}{0.5em}}

\newcommand{\po}{{\mathbb P}}
\newcommand{\epsi}{{\varepsilon}}
\newcommand{\ra}{{\mathcal R}(G)}
\begin{document}

\title{Statistical properties of stochastic 2D Navier-Stokes equations
  from linear models}
\author{Hakima Bessaih\footnote{University of Wyoming, 
Department of Mathematics, Dept. 3036, 1000
East University Avenue, Laramie WY 82071, United States, bessaih@uwyo.edu}$\;$
 \& Benedetta Ferrario\footnote{Universit\`a di Pavia, 
Dipartimento di Matematica, via
  Ferrata 5, 27100 Pavia, Italy,
benedetta.ferrario@unipv.it}}

\date{\today}
\maketitle

\begin{abstract} 
A new approach to the old-standing problem of the
  anomaly of the scaling exponents of nonlinear models of turbulence 
has been proposed and tested through numerical simulations.
This is achieved by constructing, for any given nonlinear model, 
a linear model of passive advection of an auxiliary field whose anomalous 
scaling exponents are the same as the scaling exponents of the 
nonlinear problem. 

In this paper, we investigate this conjecture for the
2D  Navier-Stokes equations driven by an additive noise.
 In order to check this conjecture, we
 analyze the coupled system Navier-Stokes/linear advection system in
 the unknowns $(u,w)$.
 We introduce a parameter $\lambda$ which
gives a system $(u^\lambda,w^\lambda)$;
this system is studied for any $\lambda$
proving its well posedness and the uniqueness of its invariant
measure $\mu^\lambda$. 

The key point is that for any $\lambda \neq 0$ 
the fields $u^\lambda$ and $w^\lambda$  have the same scaling
exponents, by assuming  universality of the scaling exponents to the
force. In order to prove the same for the original fields $u$ and $w$, 
we investigate the limit as $\lambda \to 0$, proving that
$\mu^\lambda$ weakly converges to $\mu^0$, where $\mu^0$ is the only 
invariant measure for the joint system for $(u,w)$ when $\lambda=0$.

\end{abstract}
\section{Introduction} 

We consider the stochastic Navier Stokes equations 
\begin{equation}\label{NSE}
\frac{\partial u}{\partial t}- \nu \Delta u + (u \cdot\nabla)u 
 =-\nabla p_{1}+\frac{\partial \beta_{1}}{\partial t},\qquad
\nabla \cdot u=0
\end{equation}
describing the motion of a fluid 
in a bounded domain $\mathcal D$ of $\mathbb{R}^{d}$ ($d=2,3$). 
Here,  $u=u(t,\xi)$ is the velocity vector field,
$p=p(t,\xi)$ is the scalar pressure field and  $\nu>0$ is the
viscosity. $\beta_{1}=\beta_1(t,\xi)$ is a Gaussian random field white
in time, subject to the restrictions imposed below on the space correlation.
The velocity field $u$ is subject to some  boundary condition 
and  initial condition.

Recently (see \cite{ABBPT} and the references therein),
it has been proposed that the Navier-Stokes
equation  and a relevant linear advection model  
\begin{equation}\label{Passive}
\frac{\partial w}{\partial t}- \nu \Delta w + (u \cdot\nabla)w=
-\nabla p_{2}+ \frac{\partial \beta_{2}}{\partial t}, \qquad \nabla \cdot w=0
\end{equation}
have the same
scaling exponents of their structure functions, even if their
statistics are different (i.e., in the ''jargon'' of stochastic PDE's  
their invariant measures are different). 
Here $u$ is solution of equation \eqref{NSE}; the noises
 $\beta_1$ and $\beta_2$ are 
independent and  identically distributed.
 
We recall the definition of structure function (see \cite{Fri}): it
is an ensemble average of power of velocity differences across a lenght scale.
If the velocity field $u$ is stationary in time, homogeneous and isotropic
in space,
 the longitudinal structure function is
\begin{equation}\label{s-p}
 S^p_u(l)= \Big\langle \{[u(t,l\hat \xi)-u(t,0)]\cdot \hat \xi \}^p \Big\rangle
\end{equation}
with $\xi,\xi+l\hat\xi \in D$, $p \in \mathbb N$
and $l \in \mathbb R$; $\hat\xi$ is a versor
and we take the
scalar product of vectors in $\mathbb R^d$.

For a turbulent field it is important to know how $S^p_u(l)$ 
depends on $l$ for small $l$.
The scaling exponents $\zeta_p$ are defined by
\begin{equation}
 S^p_u(l)\propto l^{\zeta_p}
\end{equation}
and are assumed to be universal in the limit $\nu \to 0$ 
when $l$ lies in the inertial range $\eta_0\nu^{3/4}\le l\le l_0$
(see e.g. \cite{Fri});
an equivalent definition is given by
\begin{equation}\label{espo-p}
 \zeta_p=\lim_{l \to 0} \frac {\log S^p_u(l)}{\log l}.
\end{equation}
Notice that $l \to 0$ implies also $\nu \to 0$.

\cite{ABBPT} provides numerical evidence that both the 2D
Navier-Stokes equations  and the Sabra shell model 
have the same scaling exponents of the corresponding linear advection
models. If this statement
were  true, then it would  allow to reduce the
establishement of the scaling exponents for the Navier-Stokes equations
\eqref{NSE} with $d=2$ 
to the easier problem of the establishement of the scaling exponents for the linear advection model
\eqref{Passive}.
As far as the 
linear advection problem is concerned, its 
scaling exponents show  anomalous behavior
(see, among the others, \cite{ABBPT,procaccia2001,GK} and the references
therein).

It is then of interest to understand rigorously the properties
of the joint system
\begin{equation}\label{Joint}
\left\{\begin{array}{lr}
 \dfrac{\partial u}{\partial t}- \nu \Delta u + (u \cdot\nabla)u 
 =-\nabla p_{1}+ \dfrac{\partial \beta_{1}}{\partial t},\\
  \dfrac{\partial w}{\partial t}- \nu \Delta w + (u \cdot\nabla)w 
 =-\nabla p_{2}+ \dfrac{\partial \beta_{2}}{\partial t}
 \end{array}\right.
\end{equation}
with the divergence free condition and the boundary condition.
If it has a unique invariant measure $\mu$, then 
the ensemble averages are computed  with respect to this measure, i.e. for
the linear advection model \eqref{Passive}
\begin{equation}\label{s-p-w}
 S^p_w(l)=\iint   \{[w(l\hat \xi)-w(0)]\cdot \hat \xi \}^p \mu(du,dw)
\end{equation}
and for the Navier-Stokes equation \eqref{NSE}
\begin{equation}\label{S-con-m}
 S^p_u(l)=\iint   \{[u(l\hat \xi)-u(0)]\cdot \hat \xi \}^p \mu(du,dw)
 \equiv \int \{[u(l\hat \xi)-u(0)]\cdot \hat \xi \}^p m(du)
\end{equation}
being $m(du)=\int \mu(du,dw)$ the unique invariant measure for \eqref{NSE}.

The analysis of  system \eqref{Joint} has been performed by adding two
terms (see \cite{BLPTiti} and the references therein): 
\begin{equation}\label{Joint_lambda}
\left\{\begin{array}{lr}
 \dfrac{\partial u}{\partial t}- \nu \Delta u + (u \cdot\nabla)u 
  +\lambda(w \cdot\nabla)u 
  =-\nabla p_{1}+ \dfrac{\partial \beta_{1}}{\partial t}
\\
  \dfrac{\partial w}{\partial t}- \nu \Delta w + (u \cdot\nabla)w
  +\lambda(w \cdot\nabla)w =-\nabla p_{2}+\dfrac{\partial \beta_{2}}{\partial t}
 \end{array}\right.
\end{equation}
where $\lambda\in \mathbb{R}$ is a parameter. We denote by 
$(u^\lambda, w^\lambda)$ its solution.
For $\lambda=0$ we recover \eqref{Joint} and for $\lambda=1$ system 
\eqref{Joint_lambda} is symmetric.

On the other hand, for any $\lambda\neq 0$, 
system \eqref{Joint_lambda} enjoys the following property:
setting $v^\lambda=\lambda w^\lambda$ and multiplying the second
equation by $\lambda$, we have a perfectly symmetric system for the
pair $(u^\lambda,v^\lambda) $, except for the force and initial
conditions:
\begin{equation}\label{coppia-u-v}
\left\{\begin{array}{lr}
 \dfrac{\partial u}{\partial t}- \nu \Delta u
  + (u \cdot\nabla)u +(v \cdot\nabla)u 
 =-\nabla p_1+\dfrac{\partial \beta_{1}}{\partial t}
\\
   \dfrac{\partial v}{\partial t}- \nu \Delta v +(u\cdot\nabla)v 
  +(v\cdot\nabla)v =-\lambda \nabla p_{2}
  +\lambda \dfrac{\partial \beta_{2}}{\partial t}
 \end{array}\right.
 \end{equation}
with $u^\lambda(0)=u_0$ and $v^\lambda(0)=\lambda w_0$.
Thus, {\it assuming}\footnote{The linear advection equation \eqref{Passive} has
scaling exponents universal to the forcing, when
the forces act only on large scales (this has been investigated by
physicists, see e.g. \cite{procaccia2001,procaccia2002,FGV}). 
Some numerical evidence
shows the same for the nonlinear equation \eqref{NSE} (see
\cite{ABBPT}).}
the universality of the scaling exponents to the
force, it follows that $u^\lambda$ and $v^\lambda=\lambda w^\lambda$ 
have the same scaling exponents for any $\lambda\neq 0$.
Then the same consideration holds for the couple $(u^\lambda,w^\lambda)$ for any
$\lambda\neq 0$; indeed, 
$S^p_{\lambda w^\lambda}(l)=\lambda^p S^p_{w^\lambda}(l)$ and
\eqref{espo-p} gives the same $\zeta_p$ for $v^\lambda$ and $w^\lambda$.
The crucial point  is to see if this holds also in
the limit as $\lambda \to 0$. 

To this end, we shall investigate when there exists a unique invariant
measure $\mu^\lambda$ 
for \eqref{Joint_lambda} (for any $\lambda\in \mathbb R$) and 
we shall prove that there exists a subsequence $\mu^{\lambda_n}$ which
converges to the unique invariant 
measure $\mu$ for system \eqref{Joint}, as $\lambda_n \to 0$. 
We recall that the
 continuous dependence of the solutions to system \eqref{Joint_lambda}
as $\lambda \to 0$ has been 
investigated rigorously in the context of certain nonlinear
phenomenological shell model (see \cite{BLPTiti}  for the deterministic
case and \cite{BFT} for the stochastic case). Moreover, \cite{BFT}
considers the asymptotic dynamics (for large time) on the attractor,
proving the continuous  dependence on $\lambda$ of the attractor.
Here we address the continuous dependence of the invariant measure with respect to
the parameter $\lambda$ for the stochastic 2D Navier-Stokes equations.
We shall exploit the fact that system \eqref{Joint_lambda}
enjoys the same properties as the stochastic 2D Navier-Stokes
equations, in order to prove existence and uniqueness of an invariant
measure $\mu^\lambda$.

The paper is organized as follows.
Section \ref{s2} is devoted to introducing the functional setting. 
It is split into subsections. We first introduce the functional spaces 
and operators. Then, the stochastic external forces and assumptions 
on the covariances of the noise are defined. Some properties of the 
Ornstein-Uhlenbeck process are introduced. 
These auxiliary and known results will be used in  Section 3
to prove the 
well posedness of \eqref{Joint_lambda} and the uniqueness of 
its invariant measure for any $\lambda\in\mathbb R$;
the proofs are given 
 along the same lines as for the Navier-Stokes equations \eqref{NSE}. 
In Section 4, we prove our main result, i.e. the continuous dependence of 
the invariant measure  $\mu^{\lambda}$ with respect to the parameter 
$\lambda$ is proved.  We first prove that the unique 
invariant measure $\mu^{\lambda}$ is tight; 
hence it  has a limit $m$ when $\lambda\to 0$. Now, in order to prove 
that the limit $m$ is an invariant measure, we consider the stationary 
solutions of \eqref{Joint_lambda} whose time marginals are
$\mu^{\lambda}$. 
Thanks to uniform estimates computed on these stationary solutions, 
we prove their convergence when $\lambda\to 0$ and that their limit 
are stationary and  their marginals are given by $m$. 
As a consequence, we conclude that $m$ is the unique invariant
measure for  \eqref{Joint_lambda} when $\lambda=0$.

%____________________________________________________________________
\section{Notations and hypothesis}\label{s2}
\subsection{Functional setting}
We define the functional setting to study the Navier-Stokes
equations. From now on, the spatial domain is the square 
$\mathcal D=[0,2\pi]^2$ with periodic boundary conditions.  

As usual, in the periodic case we assume that the spatial 
mean of the vectors we are dealing with is zero. 
This gives a simplification in the mathematical treatment, 
but it does not prevent to consider non zero mean value vectors. 
Actually, if we can analyze the problem for zero mean vectors 
then the problem without this assumption can be dealt with in 
a similar way (see \cite{Tper}).

Let $H^s$ be the space of divergence free and periodic vector fields with mean zero that 
 belong to  the  Sobolev space $[H^s(\mathcal D)]^2$.
\\
Let 
$\Pi$ be the orthogonal projector in $[L^{2}({\mathcal  D})]^{2}$ onto $H^0$;
then the  Stokes operator is defined as 
$$ 
 Au=-\Pi \Delta u, \qquad \forall u \in  D(A)=H^2.
$$
It is a closed positive unbounded self-adjoint operator in $H$ with
the inverse 
$A^{-1}$ which is a 
self-adjoint compact operator in $H^0$; by the classical
spectral theorems there exists 
a sequence $\{\gamma_j\}_{j=1}^{\infty}$ 
of eigenvalues of the Stokes operator  
with $0<\gamma_{1}\le \gamma_{2}\le \dots$, 
corresponding to the eigenvectors $e_{j} \in D(A)$; $\{e_j\}_{j\in \mathbb N}$ 
form an orthonormal basis in $H^0$.
We have that $\gamma_j$ 
behaves like $j$ for $j\to \infty$.
The Stokes operator generates an analytic semigroup $e^{-\nu tA}$ in $H^0$ and 
for each $s>0$ there exists a 
constant $M>0$ (depending on $s$ and $\nu$) such that
\[
\|e^{-\nu tA}u\|_{H^s}\le \frac M{t^{s/2}} \|u\|_{H^0} , \qquad u \in H^0
\]
for all $t>0$.

We can define the fractional powers $A^p$ ($p>0$) as
linear unbounded operators in $H^0$ with 
\[
A^pe_j=\gamma_j^p e_j, \qquad D(A^p)=H^{2p}.
\]
Therefore we can characterize the Hilbert spaces $H^s$ as
\[
H^s=\{u=\sum_{j=1}^\infty u_j e_j: \sum_{j=1}^\infty \gamma_j^s u_j^2<\infty\}
\]
and we set
\[
\|u\|_{H^s}^2= \sum_{j=1}^\infty \gamma_j^s u_j^2.
\]
Moreover,  $H^{s_1}$ is densely and compactly embedded in 
$H^{s_2}$ for $s_1>s_2$.
Finally, $H^{-s}$ denotes the dual of $H^s$, with duality bracket
$\langle\cdot,\cdot \rangle$. For $s=0$ this is the scalar
product $\langle u,v\rangle=\sum_j u_j v_j$ in $H^0$.

Let $b(\cdot,\cdot,\cdot)$ be the  trilinear form defined as
$$
 b(u,v,z)=\int_{\mathcal D}([u(\xi)\cdot\nabla] v(\xi))\cdot z(\xi)\, d\xi .
$$
It is well known that there exists a continuous bilinear operator
$B(\cdot,\cdot): H^1\times H^1\longrightarrow H^{-1}$ such that
$\langle  B(u,v),z\rangle =b(u,v,z)$ for all $z\in H^1$.
By the incompressibility condition, 
we have 
\begin{equation} \label{incompress}
 \langle B(u,v),z\rangle=- \langle B(u,z),v\rangle
 \quad   \mbox{\rm and}\quad    \langle B(u,v),v\rangle =0
\end{equation}
for $u,v, z\in H^1$. 
Furthermore, there exist  constants $C$ and $C_r$ such that 
\begin{equation}\label{bilinear_estimate1}
\begin{cases}
&
 \|B(u,v)\|_{H^{-1}}\le C \|u\|_{H^{\frac 12}}\|v\|_{H^{\frac
     12 }}\\
     &
 \|B(u,v)\|_{H^{0}}\le C \|u\|_{H^{\frac 12}}\|v\|_{H^{\frac
     32 }}\\
\text{ for } r \ge 2&
\|B(u,v)\|_{H^{-1+r}}\le C_r \|u\|_{H^r} \|v\|_{H^r}
\end{cases}
\end{equation}
A more refined inequality holds: for  $r>2$
\[
\|B(u,v)\|_{H^{-1+r}}\le C \|u\|_{H^{-1+r}} \|\nabla v\|_{H^{-1+r}}
\le C \|u\|_{H^{-1+r}} \|v\|_{H^r}
\]
since $H^{-1+r}$ is a multiplicative algebra for $r>2$.
Since 
\[\begin{split}
\|B(u,v)\|_{H^1}&=\|(u \cdot \nabla)v\|_{H^1}\le \|Du\|_{L^4}\|\nabla
v\|_{L^4}+ \|u\|_{L^\infty}\|\nabla v\|_{H^1}
\\&\le
C \|u \|_{H^{\frac 32}}\|v\|_{H^{\frac 32}}+C \|u\|_{H^{\frac 32}}\|v\|_{H^2}
\\&
\le C \|u\|_{H^{\frac 32}}\|v\|_{H^2}
\end{split}\]
these two latter inequalities are summarized in the last line 
of \eqref{bilinear_estimate1}.

Finally, let us point out that we will use the same symbol $C$ for
different constants, if they do not play an important role.

%%%%%%%%%%%%%%%%%%
\subsection{Stochastic driving force}

As far as the stochastic forcing terms are concerned,
we refer to \cite{DaPZ} for the basic results. 
Here we recall the main definitions and properties.

We introduce two independent $H^0$-cylindrical Wiener processes
$\beta_1$ and $\beta_2$ as follows:
\begin{equation}\label{W-n}
\beta_{i}(t)=\sum_{j=1}^{\infty} \omega^{(i)}_j(t) e_j,
\end{equation}
with  $(\omega^{(i)}_j)$   mutually independent standard (scalar)
Wiener processes 
defined   on a filtered  complete probability space 
$(\Omega, \mathcal{F}, (\mathcal{F}_t)_{t\ge 0}, \mathbb{P})$.\\
Let $G$ be a  Hilbert-Schmidt operator in $H^0$ and denote by 
$\ra$ its range and by $\|G\|_{HS}$ its
Hilbert-Schmidt norm ($\|G\|_{HS}^2=\sum_j \|Ge_j\|^2_{H^0}$).
If we assume that 
\[
\ra \subseteq H^p
\]
for some $p\ge 0$, then 
each process $G\beta_i$ takes values in $C(\mathbb R_+;H^s)$ if
\begin{equation}\label{ccw}
p>s+1 .
\end{equation}
Indeed this is 
equivalent to  $A^{\frac s2}G\beta_i$ taking values in $C(\mathbb R_+;H^0)$.
Notice that $A^{\frac s2}G= A^{\frac{s-p}2} (A^{\frac p2}G)$;
therefore, by assuming that  $A^{\frac p2}G$
is a bounded operator in $H^0$, we get that $A^{\frac s2}G$ is a
Hilbert-Schmidt operator if $A^{\frac{s-p}2}$ is so. This happens
under assumption \eqref{ccw}, since 
$\|A^{\frac{s-p}2}\|^2_{HS}=\sum_j \|A^{\frac{s-p}2}e_j\|^2_{H^0}=
\sum_j \gamma_j^{s-p}\sim \sum_j j^{s-p}$ and the latter series is
convergent iif \eqref{ccw} is fulfilled.

From now on, we assume that 
\[
\exists \ \epsi>0: \; \ra \subseteq H^{1+\epsi}.
\]
This implies that $G\beta_i \in C(\mathbb R_+;H^0)$, $\po$-a.s..

Projecting the equations of \eqref{Joint_lambda}  onto $H^0$, 
we get rid of the pressure terms and the following
abstract formulation is obtained
\begin{equation}
\label{Joint_abstract}\left\{
\begin{array}
[c]{ll}%
du^{\lambda}(t)+[\nu Au^{\lambda}(t)+B(u^{\lambda}(t),u^{\lambda}(t))+\lambda
B(w^{\lambda}(t),u^{\lambda}(t))]dt= G d\beta_{1}(t), & \quad t>0\\
dw^{\lambda}(t)+[\nu Aw^{\lambda}(t)+B(u^{\lambda}(t),w^{\lambda}(t))+\lambda
B(w^{\lambda}(t),w^{\lambda}(t))]dt= G d\beta_{2}(t), & \quad t>0\\
u^{\lambda}(0)=u_{0} & \\
w^{\lambda}(0)=w_{0} &
\end{array}
\right. 
\end{equation}
with initial conditions $u_{0},w_{0}\in H^0$.
\\
The Cauchy problem is studied on any finite time interval $[0,T]$.

%%%%%%%%%%%%%%%%%%%%%%%%%%%%%%%%%
\subsection{Compact form of the system}

Now we write this system in a compact form.
Define $\widetilde{H}^s=H^s\times H^s$.
If $x=(x_{1},x_{2})\in\widetilde{H}^0$ and
$y=(y_{1},y_{2})\in\widetilde{H^0}$, we define the scalar product in
$\widetilde{H}^0$ as
\[
<x,y>=<x_{1},y_{1}>+<x_{2},y_{2}>
\]
and the norms in $\widetilde{H}^s$ as
\[
\|x\|_{\widetilde{H}^s}^{2}=\|x_{1}\|_{H^s}^{2}+\|x_{2}\|_{H^s}^{2},\quad x=(x_{1}%
,x_{2})\in\widetilde{H}^s .
\]

We define the linear operator $\widetilde{A}:\widetilde H^2
 \subset\widetilde{H}^0\rightarrow\widetilde{H}^0$ as 
$\tilde A x =\big(Ax_1, Ax_2\big)$. As a consequence, $\tilde A$ is nonnegative and selfadjoint. 

For every $\lambda\in \mathbb{R}$, we define the bilinear operator
$\widetilde{B}^{\lambda}$  as
\begin{equation}\label{bt}
 \widetilde{B}^{\lambda}\left(  x,y\right)  
 =\Big(  B(x_{1},y_{1})  +\lambda B(x_{2},y_{1}),B(x_{1},y_{2})
   +\lambda B(x_{2},y_{2})  \Big),
\end{equation}
By \eqref{incompress}, \eqref{bilinear_estimate1} and \eqref{bt}, we have
\begin{lemma}\label{lemma su B tilde}
For any $x,y,z \in \tilde H^1$ we have
\[
 \langle \widetilde{B}^{\lambda}(x,y),z\rangle=
 -\langle \widetilde{B}^{\lambda}(x,z),y\rangle,
\qquad
 \langle \widetilde{B}^{\lambda}(x,y),y\rangle=0.
\]
Moreover\\
i) there is a constant $C>0$ 
such that for any  $\lambda\in \mathbb R$ 
\[
\|\widetilde{B}^{\lambda}(x,y)\|_{\widetilde{H}^{-1}}
\le C(1+|\lambda|) 
\|x\|_{\tilde H^{\frac 12}} \|y\|_{\tilde H^{\frac 12}};
\]
\\ii) 
there is a constant $C>0$ 
such that for any  $\lambda\in \mathbb R$ 
\[
\|\widetilde{B}^{\lambda}(x,y)\|_{\widetilde{H}^{0}}
\le C (1+|\lambda|)
\|x\|_{\tilde H^{\frac 12}} \|y\|_{\tilde H^{\frac 32}};
\]
\\iii) for any $r\ge 2$
there is a constant $C_{r}>0$ 
such that for any $\lambda\in \mathbb R$
\[
\|\widetilde{B}^{\lambda}(x,y)\|_{\widetilde{H}^{-1+r}}
\le C_{r} (1+|\lambda|)
\|x\|_{\tilde H^{r}} \|y\|_{\tilde H^{r}}.
\]
\end{lemma}
Actually the latter estimate comes from
\begin{equation}\label{perdimo1}
\|\widetilde{B}^{\lambda}(x,y)\|_{\widetilde{H}^{1}}
\le C (1+|\lambda|)
\|x\|_{\tilde H^{\frac 32}} \|y\|_{\tilde H^{2}}
\end{equation}
and, for $r \ge 3, r \in \mathbb N$ 
\begin{equation}\label{perdimo2}
\|\widetilde{B}^{\lambda}(x,y)\|_{\widetilde{H}^{-1+r}}
\le C_{r} (1+|\lambda|)
\|x\|_{\tilde H^{r-1}} \|y\|_{\tilde H^{r}}.
\end{equation}

Hence, we write \eqref{Joint_abstract} in more compact form as
\begin{equation}\label{NSE_lambda}
 \left\{\begin{array}[c]{ll}
 d\widetilde{u}^{\lambda}(t)+[  \nu\widetilde{A}\widetilde{u}^{\lambda}(t)
 +\widetilde{B}^{\lambda}\left(  \widetilde{u}^{\lambda}(t),\widetilde
 {u}^{\lambda}(t)\right)  ]  dt  = \tilde G
 d\widetilde{\beta}(t),\qquad t>0\\
 \widetilde{u}^{\lambda}(0)  =\widetilde x
\end{array}
\right.
\end{equation}
where $\widetilde{u}^{\lambda} =(u^{\lambda},w^{\lambda})$ and
$\tilde G \tilde \beta=\big(G \beta_1, G\beta_2\big)$.

The following definition of solution for \eqref{NSE_lambda}
is given 
\begin{definition}
A stochastic process $\tilde u^\lambda$ is a generalized solution
in $[0,T]$ of system \eqref{NSE_lambda} if
\[
\tilde u^\lambda \in C([0,T];\tilde H^0) \cap L^2(0,T;\tilde H^1)
\]
$\po$-a.s.  and equation \eqref{NSE_lambda}
 is satisfied $\mathbb P$-a.s. in the integral
sense
\[
\langle\tilde u^\lambda(t),\phi\rangle
+\nu  \int_0^t 
\langle\tilde A^{\frac 12}\tilde u^\lambda(s),\tilde A^{\frac
  12}\phi\rangle ds 
+
\int_0^t \langle\tilde B^\lambda(\tilde u^\lambda(s),\tilde u^\lambda(s)),
\phi\rangle \;ds =
\]
\[
= \langle \tilde{x}, \phi\rangle +\langle \tilde G \tilde \beta(t),\phi\rangle
\]
for all $t \in [0,T]$ and all $\phi \in \tilde H^1$.
\end{definition}

This definition involves stronger regularity than in 
\cite{Ferrario03}-\cite{Fladissipative}, but in those papers one
looked for the minimal assumption on the noise to define a
solution. But here we are interested in more regular solutions; for
this reason we shall assume that $G$ is an Hilbert-Schmidt operator in
$H^0$  and so we avoid some technicalities of
\cite{Ferrario03}-\cite{Fladissipative}, required to deal with 
solutions of low space  regularity.

%%%%%%%%%%%%%%%%%%
\subsection{Ornstein-Uhlenbeck Process}\label{sdf}
Let us define the Ornstein-Uhlenbeck equation
\begin{equation}\label{eqOU}
dz_i(t)+ \nu A z_i(t)\ dt = G d \beta_i(t).
\end{equation}
We consider its mild solution
\begin{equation}\label{procOU}
z_i(t)=e^{-\nu At}z_i(0) +\int_{0}^t e^{-\nu A(t-s)}Gd\beta_i(s)
\end{equation}
whose regularity is given in ths following proposition (see
\cite{dpkz}).
\begin{proposition}\label{p:ou}
Let $G:H^0\to H^0$ be a linear bounded operator such that 
 $\ra\subseteq H^{\alpha+\epsi}$ for some $\alpha, \epsi>0$.
Then, given any initial value  $z_i(0) \in H^\alpha$, there exists a
version of the Ornstein-Uhlenbeck process \eqref{procOU}
whose   paths are  in $C([0,\infty[;H^\alpha)$, $\po$-a.s..
\end{proposition}

Therefore, from now on our  assumption on the covariance of the noises is given by
\begin{description}
\item[ {$[$H0$]$} ]
$G:H^0\to H^0$ is a linear bounded operator such that 
 $\ra\subseteq H^{1+\epsi}$ for some $\epsi>0$
\end{description}
This implies that $G$ is an Hilbert-Schmidt operator in $H^0$.
 
Notice here that the mild solution of the Ornstein-Uhlenbeck equation is also a weak solution (in the PDE sense).
 We refer to \cite{DaPZ} (Theorem 5.4). 

%%%%%%%%%%%%%%%%%%%%%%
\section{Ergodicity}\label{s3}
In this section we deal with fixed $\lambda$, arbitrary chosen in
$\mathbb R$.
For the evolution problem \eqref{NSE_lambda},
we will prove first that if the initial data 
$\tilde x\in \tilde H^0$ and {\bf H0}
is fulfilled, then there exists a unique generalized solution, which
is a strong solution in the probabilistic sense.
This unique solution is a time-homogeneous Markov process in $\tilde H^0$.
Then,  we  define the  transition functions 
$P^\lambda(t,\tilde x,\Gamma) := 
\po \{\tilde u^\lambda(t;\tilde u^\lambda(0)=\tilde x) \in \Gamma \}$
for any  $t \ge 0$, $\tilde x \in \tilde H^0$ and 
Borel subset $\Gamma$ of $\tilde H^0$;
the  associated Markov semigroup is
\[
(P^\lambda_t \psi)(\tilde x)
=
\mathbb E [\psi(u^\lambda(t;\tilde u^\lambda(0)=\tilde x))].
\]
For each $t$, $P^\lambda_t$  maps 
Borelian bounded functions $B_b(\tilde H^0)$ in themselves; it 
is said to be  Feller if it maps 
continuous bounded functions in themselves, i.e.
$P^\lambda_t:C_b(\tilde H^0) \rightarrow C_b(\tilde H^0)$. 
Going further, invariant measures will be investigated,
 that is measures $\mu^\lambda$ such that
\[
\int \psi \ d\mu^\lambda = \int P^\lambda_t\psi \ d\mu^\lambda
\] 
for any $t>0$ and $\psi \in C_b(\tilde H^0)$. 
We will prove that there exists a unique invariant measure for
equation \eqref{NSE_lambda}.

%%%%%%%%%%%%%%%%%%%%%%%%%%%%
\subsection{Well posedness}
In this section we prove existence, uniqueness and continuous
dependence of the solution on the intial data for system 
\eqref{Joint_lambda}.
The result is a classical one, 
since we are in a two dimensional spatial domain.
We shall first prove existence of a martingale 
solution.
Then pathwise uniqueness and continous
dependence of the solution on the intial data. 
Pathwise uniqueness implies that we have actually a strong solution in
the probabilistic sense.

Our procedure is to work pathwise, as in \cite{Fladissipative} and \cite{Ferrario03},
whereas \cite{Chow} and \cite{VF} (Ch X) look for mean square
estimates. In this way we get better regularity for the paths, 
i.e. $\tilde u^\lambda \in  C([0,T];\widetilde{H}^0)$ and not only 
$\tilde u^\lambda \in  L^\infty(0,T;\widetilde{H}^0)$. 
\begin{proposition} \label {esistH}
Assume {\bf H0}.
\\
Then, for any time interval $[0,T]$ and  any initial condition 
$x \in \widetilde{H}^0$
there exists a unique strong generalized solution $\tilde u^\lambda$ of
equation \eqref{NSE_lambda} over $[0,T]$ with the initial condition
$\tilde u^\lambda(0)=\tilde{x}$, satisfying $\po$-a.s.
$$ 
 \tilde u^\lambda \in C([0,T];\widetilde{H}^0) \cap  
         L^2(0,T;\widetilde{H}^1). 
$$
The process $\tilde u^\lambda$ is a Markov and Feller process in
$\widetilde{H}^0$.
\end{proposition}
\proof
The proof relies on a priori estimates for the  Ornstein-Uhlenbeck
process $\tilde z$, solving the linear equation
\[
d\tilde z(t)+\nu \tilde A \tilde z(t)\ dt 
=\tilde G d\tilde w(t); \qquad \tilde z(0)=0
\]
(so $\tilde z=(z_1, z_2)$ with $z_i \in C([0,\infty[;H^1))$ a.s., 
see \eqref{eqOU})
and the auxiliary process
$\tilde v^\lambda=\tilde u^\lambda-\tilde z$, 
solving the deterministic equation 
\[
\frac {d\tilde v^\lambda}{dt}(t)+\nu \tilde A \tilde v^\lambda (t)
+\tilde B^\lambda(\tilde v^\lambda(t)+\tilde z(t), \tilde v^\lambda(t)+\tilde
z(t))=0;
\qquad \tilde v^\lambda(0)=\tilde{x}
\]
with random coefficients. 
The equation for $\tilde v^\lambda$ is equivalent to
\begin{multline}\label{eq-v}
\frac {d\tilde v^\lambda}{dt}(t)+\nu \tilde A \tilde v^\lambda (t)
+\tilde B^\lambda(\tilde v^\lambda(t),\tilde v^\lambda(t))
+\tilde B^\lambda(\tilde v^\lambda(t),\tilde z(t))
+\tilde B^\lambda(\tilde z(t), \tilde v^\lambda(t))
\\
=-\tilde B^\lambda(\tilde z(t), \tilde z(t))
\end{multline}
Local existence is easy to prove, since the nonlinearity is locally
Lipschitz. Therefore it is enough to get apriori estimates to show 
that this solution exists on the whole time interval $[0,T]$.
Actually, one should first work with the finite-dimensional
(Galerkin) approximation and then pass to the limit; we omit this
details since they are classical (see, e.g.,  \cite{Temam}, \cite{Tper}).

We take the $\tilde H^0$-scalar product of equation \eqref{eq-v}
with $\tilde v^\lambda$ and use that 
$\tilde H^{1/2}\subset [L^4(\mathcal D)]^2$ and by interpolation
$\|v\|_{\tilde H^{1/2}}\le C \|v\|_{\tilde H^0}^{1/2} \|v\|_{\tilde H^1}^{1/2}$.
Therefore, by taking into account Lemma \ref{lemma su B tilde}
we get
\[\begin{split}
\tfrac 12 \tfrac {d}{dt}\|\tilde v^\lambda(t)\|_{\tilde H^0}^2
&+\nu \|\tilde v^\lambda(t)\|^2_{\tilde H^1}
=
-\langle \tilde B^\lambda(\tilde v^\lambda(t),\tilde z(t)),
\tilde v^\lambda(t) \rangle
%\\&\qquad\qquad\qquad\qquad\qquad
-\langle \tilde B^\lambda(\tilde z(t),\tilde z(t)),
\tilde v^\lambda(t) \rangle
\\&
\le (1+|\lambda|)\|\tilde v^\lambda(t)\|_{L^4}^2 \|\tilde z(t)\|_{\tilde H^1}+
(1+|\lambda|)\|\tilde z(t)\|_{L^4} \|\tilde z(t)\|_{\tilde H^1}
\|\tilde v^\lambda(t)\|_{L^4}
\\
&\le
C (1+|\lambda|)\|\tilde v^\lambda(t)\|_{\tilde H^0} \|\tilde v^\lambda(t)\|_{\tilde H^1} 
\|\tilde z(t)\|_{\tilde H^1}+
C (1+|\lambda|)\|\tilde v^\lambda(t)\|_{\tilde H^1} \|\tilde z(t)\|_{\tilde H^1}^2
\\& \le
\frac \nu 2 \|\tilde v^\lambda(t)\|_{\tilde H^1}^2
+ C_\nu (1+\lambda^2)\|\tilde z(t)\|_{\tilde H^1}^2 
\|\tilde v^\lambda(t)\|_{\tilde H^0}^2
+C_\nu (1+\lambda^2) \|\tilde z(t)\|_{\tilde H^1}^4
\end{split}
\]
where we have used Young inequality in the latter step.
Thus
\begin{equation}\label{stima1}
\frac {d}{dt}\|\tilde v^\lambda(t)\|_{\tilde H^0}^2
+\nu \|\tilde v^\lambda(t)\|^2_{\tilde H^1}
\le 
 C\|\tilde z(t)\|_{\tilde H^1}^2
\|\tilde v^\lambda(t)\|_{\tilde H^0}^2
+C \|\tilde z(t)\|_{\tilde H^1}^4
\end{equation}
where the constant $C$ depends on $\nu$ and $\lambda$ and is uniformly
 bounded for 
$\lambda$ in a bounded set.
Gronwall lemma applied to 
\[
\frac {d}{dt}\|\tilde v^\lambda(t)\|_{\tilde H^0}^2\le 
C\|\tilde z(t)\|_{\tilde H^1}^2
\|\tilde v^\lambda(t)\|_{\tilde H^0}^2
+C \|\tilde z(t)\|_{\tilde H^1}^4
\]
gives
\[
\sup_{0\le t \le T}\|\tilde v^\lambda(t)\|_{\tilde H^0}^2\le
\|\tilde{x}\|_{\tilde H^0}^2
e^{\int_{0}^T C\|\tilde z(t)\|_{\tilde H^1}^2 dt}
+C \int_{0}^T 
e^{\int_{0}^s C\|\tilde z(r)\|_{\tilde H^1}^2 dr }
\|\tilde z(s)\|_{\tilde H^1}^4 ds 
\]
where the r.h.s. is finite, since {\bf H0} grants that 
the paths of $\tilde z$ are in 
$C([0,T];\tilde H^1)$, $\po$-a.s..
Integrating in time \eqref{stima1}
we get
\[
\nu \int_{0}^T \|\tilde v^\lambda(t)\|^2_{\tilde H^1}dt\le 
\|\tilde{x}\|_{\tilde H^0}^2+
C\left( \sup_{0\le t \le T}\|\tilde v^\lambda(t)\|_{\tilde H^0}^2 \right)
\int_{0}^T \|\tilde z(t)\|_{\tilde H^1}^2 dt+C 
\int_{0}^T \|\tilde z(t)\|_{\tilde H^1}^4 dt
\]
where the r.h.s. is finite.

In this way one gets $\tilde v^\lambda\in L^\infty(0,T;\tilde H^0)\cap
L^2(0,T;\tilde H^1)$ and $\tilde u^\lambda=\tilde v^\lambda+\tilde z
\in L^\infty(0,T;\tilde H^0)\cap L^2(0,T;\tilde H^1)$. 
The continuity in time comes from the time
regularity of $\tilde v^\lambda$. First we notice that 
\[\begin{split}
\|\tilde B^\lambda(u,u)\|_{\tilde H^{-1}}=\sup_{\|\phi\|_{\tilde
    H^1}\le1}
|\langle \tilde B^\lambda(u,u), \phi \rangle | 
&=\sup_{\|\phi\|_{\tilde
    H^1}\le1}
|\langle \tilde B^\lambda(u,\phi), u \rangle | 
\\&\le C (1+|\lambda|)
\|u\|_{\tilde H^0}\|u\|_{\tilde H^1}
\end{split}\]
and therefore
$\tilde B^\lambda(\tilde u^\lambda,\tilde u^\lambda)
\in L^2(0,T;\tilde H^{-1})$.
Thus for the time derivative we get
 $\frac {d\tilde v^\lambda}{dt}(t)=-\tilde A \tilde v^\lambda (t)
-\tilde B^\lambda(\tilde u^\lambda(t),\tilde u^\lambda(t))
 \in L^2(0,T;\tilde H^{-1})$. This implies (see \cite{Temam}) that
$\tilde v^\lambda \in C([0,T];\tilde H^0)$.

Given the existence of a process $\tilde u^\lambda$ whose paths are 
in $C([0,T];\tilde H^0)\cap L^2(0,T;\tilde H^1)$, $\po$-a.s., 
we easily get pathwise uniqueness and continuous dependence on the
initial data. Let 
$\tilde u^{\lambda,(1)}$ and $\tilde u^{\lambda,(2)}$ be two
solutions of \eqref{NSE_lambda} with initial data  $\tilde x^{(1)}$ and
$\tilde x^{(2)}$  respectively.
Denote their difference by 
$U^\lambda =\tilde u^{\lambda,(1)}-\tilde u^{\lambda,(2)}$.
Since the noise is additive, the unknown $U^\lambda$ satisfies
a deterministic equation:
\[
\frac{dU^\lambda}{dt} (t)+\nu \tilde A U^\lambda(t)
+\tilde B^\lambda(\tilde u^{\lambda,(1)}(t),\tilde u^{\lambda,(1)}(t))-
\tilde B^\lambda(\tilde u^{\lambda,(2)}(t),\tilde u^{\lambda,(2)}(t))=0
\]
i.e.
\begin{equation}\label{eq:diff}
\frac{dU^\lambda}{dt}(t)+\nu \tilde A U^\lambda(t)
+\tilde B^\lambda(U^{\lambda}(t),\tilde u^{\lambda,(1)}(t))+
\tilde B^\lambda(\tilde u^{\lambda,(2)}(t),U^{\lambda}(t))=0.
\end{equation}
With estimates as before, we get
\[\begin{split}
\frac 12 \frac {d}{dt} \|U^\lambda(t)\|_{\tilde H^0}^2
+\nu  \|U^\lambda(t)\|_{\tilde H^1}^2
&=
-\langle \tilde B^\lambda(U^{\lambda}(t),\tilde u^{\lambda,(1)}(t)), 
U^\lambda(t)\rangle
\\
&\le (1+|\lambda|) \|U^{\lambda}(t)\|_{L^4}^2 
\|\tilde u^{\lambda,(1)}(t)\|_{\tilde H^1}
\\
&\le
C (1+|\lambda|) \|U^{\lambda}(t)\|_{\tilde H^0} \|U^{\lambda}(t)\|_{\tilde H^1}
\|\tilde u^{\lambda,(1)}(t)\|_{\tilde H^1}
\\&
\le \frac \nu 2 \|U^{\lambda}(t)\|_{\tilde H^1}^2+C_\nu (1+\lambda^2)
\|\tilde u^{\lambda,(1)}(t)\|_{\tilde H^1}^2 \|U^{\lambda}(t)\|_{\tilde H^0}^2
\end{split}\]
From Gronwall lemma we get
\[
\sup_{t\in [0,T]} \|U^\lambda(t)\|_{\tilde H^0}^2
\le \|U^\lambda(0)\|_{\tilde H^0}^2 
e^{C_\nu (1+\lambda^2)\int_{0}^T \|\tilde
u^{\lambda,(1)}(s)\|_{\tilde H^1}^2 ds}
\]
giving continuous dependence on the initial data. When 
$U^\lambda(0)=0$, this gives $\|U^\lambda(t)\|_{\tilde H^0}=0$ for all
$t>0$
and proves pathwise uniqueness.

The pathwise continuous dependence on the initial
data 
\[
\tilde u^{\lambda,(2)} \to \tilde u^{\lambda,(1)}\text{ in }
 C([0,T];\tilde H^0) \qquad \text{ if }
x^{(2)} \to x^{(1)} \text{ in } \tilde H^0
\]
gives the Feller property in $\tilde H^0$. Indeed, given a continuous
function $\psi: \tilde H^0\to \mathbb R$ we have 
$\po$-a.s.  that $\psi(\tilde u^{\lambda,(2)}(t))\to
\psi(\tilde u^{\lambda,(1)}(t))$  as 
$x^{(2)} \to x^{(1)}$ in  $\tilde H^0$; since $\psi$ is bounded, by dominated convergence we
get the convergence also when we take the mathematical
expectation, i.e. 
$P^\lambda_t\psi(x^{(2)})=\mathbb E [\psi(\tilde u^{\lambda,(2)}(t))]\to
\mathbb E [\psi(\tilde u^{\lambda,(1)}(t))]=P^\lambda_t\psi(x^{(1)})$.
\hfill $\Box$

\medskip
We have a regularity result.

\begin{proposition}
Assume $\ra \subseteq H^{\frac 32+\epsi}$ for some $\epsi >0$.
Then for arbitrary $\tilde{x}\in \tilde H^1$,
there exists  a unique process $\tilde u^\lambda$ solution to system 
\eqref{NSE_lambda}  over $[0,T]$ with the initial condition
$\tilde u^\lambda(0)=\tilde{x}$, such that
\[
 \tilde u^\lambda \in C([0,T];\tilde H^1)\cap L^4(0,T;\tilde H^{\frac 32}) 
\] 
$\po-$a.s.; 
it is a Markov  and Feller process in $\tilde H^1$.
\end{proposition}
\proof By assumption we have $\tilde z \in C([0,T];\tilde H^{\frac 32})$,
$\po$-a.s..
We study the regularity of $\tilde v^\lambda$.
We take the $\tilde H^0$-scalar product of equation \eqref{eq-v}
with $\tilde A \tilde v^\lambda$. Then
\[\begin{split}
\tfrac 12 \tfrac {d}{dt}\|\tilde v^\lambda(t)\|_{\tilde H^1}^2
&+\nu \|\tilde v^\lambda(t)\|^2_{\tilde H^2}
\\&=
-\langle 
\tilde B^\lambda(\tilde v^\lambda(t)+\tilde z(t),\tilde v^\lambda(t)+\tilde z(t)),
\tilde A\tilde v^\lambda(t) \rangle
\\&\le
\|\tilde B^\lambda(\tilde v^\lambda(t)+\tilde z(t),\tilde v^\lambda(t)+\tilde z(t))\|_{\tilde H^0} \|\tilde v^\lambda(t)\|_{\tilde H^2}
\end{split}
\]

We use the bilinearity, the estimates of Lemma \ref{lemma su B tilde}
and the interpolation estimates:
\[\begin{split}
\|\tilde B^\lambda(\tilde v^\lambda+\tilde z,&\tilde v^\lambda+\tilde z)\|_{\tilde H^0}
\\&
\le C (1+|\lambda|)
(\|\tilde v^\lambda\|_{\tilde H^{1/2}}\|\tilde v^\lambda\|_{\tilde H^{3/2}}+
\|\tilde v^\lambda\|_{\tilde H^{3/2}}\|\tilde z\|_{\tilde H^{3/2}}+
\|\tilde z\|^2_{\tilde H^{3/2}}) 
\\&
\le C (1+|\lambda|)
(\|\tilde v^\lambda\|_{\tilde H^{0}}^{1/2}\|\tilde v^\lambda\|_{\tilde H^{1}}
\|\tilde v^\lambda\|_{\tilde H^{2}}^{1/2}+
\|\tilde v^\lambda\|_{\tilde H^{1}}^{1/2}\|\tilde v^\lambda\|_{\tilde H^{2}}^{1/2}
\|\tilde z\|_{\tilde H^{\frac 32}}+
\|\tilde z\|^2_{\tilde H^{\frac 32}})
\end{split}
\]
So, by Young inequality
\[\begin{split}
\|\tilde B^\lambda(\tilde v^\lambda+\tilde z,&\tilde v^\lambda+\tilde z)\|_{\tilde H^0} \|\tilde v^\lambda\|_{\tilde H^2}
\\&\le
C (1+|\lambda|)
(\|\tilde v^\lambda\|_{\tilde H^0}^{1/2}\|\tilde v^\lambda\|_{\tilde H^{1}}
\|\tilde v^\lambda\|_{\tilde H^{2}}^{3/2}+
\|\tilde v^\lambda\|_{\tilde H^{1}}^{1/2}\|\tilde v^\lambda\|_{\tilde H^{2}}^{3/2}
\|\tilde z\|_{\tilde H^{\frac 32}}+
\|\tilde z\|^2_{\tilde H^{\frac 32}} \|\tilde v^\lambda\|_{\tilde H^{2}})
\\&\le
\frac \nu 2 \|\tilde v^\lambda\|_{\tilde H^{2}}^2+C_\nu(1+\lambda^2)\left[
(\|\tilde v^\lambda\|_{\tilde H^0}^{2}\|\tilde v^\lambda\|_{\tilde H^{1}}^2+
\|\tilde z\|_{\tilde H^{\frac 32}}^4)\|\tilde v^\lambda\|_{\tilde H^{1}}^2
+\|\tilde z\|_{\tilde H^{\frac 32}}^4\right]
\end{split}
\]
Therefore, with 
$g:= \|\tilde v^\lambda\|_{\tilde H^0}^{2}\|\tilde v^\lambda\|_{\tilde H^{1}}^2+
\|\tilde z\|_{\tilde H^{\frac 32}}^4\in L^1(0,T)$ by the previous
Proposition, we have
\begin{equation}\label{eq-v-H1}
\tfrac {d}{dt}\|\tilde v^\lambda(t)\|_{\tilde H^1}^2
+\nu \|\tilde v^\lambda(t)\|^2_{\tilde H^2}
\le
C_\nu(1+\lambda^2)
g(t)\|\tilde v^\lambda(t)\|_{\tilde H^{1}}^2
+C_\nu(1+\lambda^2)\|\tilde z(t)\|_{\tilde H^{\frac 32}}^4
\end{equation}
and as before we get
\[
\sup_{0\le t \le T}\|\tilde v^\lambda(t)\|_{\tilde H^{1}}^2
\le
\|\tilde{x}\|_{\tilde H^{1}}^2e^{C_\nu(1+\lambda^2)\int_0^T g(t)dt}
+C_\nu(1+\lambda^2)\int_0^T e^{C_\nu(1+\lambda^2)\int_0^T
  g(t)dt}\|\tilde z(t)\|_{\tilde H^{\frac 32}}^4dt
\]
and integrating in time \eqref{eq-v-H1} we get
\[
\nu \int_0^T\|\tilde v^\lambda(t)\|^2_{\tilde H^1}dt\le  \|\tilde{x}\|_{\tilde H^{1}}^2
+C_\nu(1+\lambda^2)\left[\sup_{0\le t \le T}\|\tilde
  v^\lambda(t)\|_{\tilde H^{1}}^2\int_0^T g(t)dt
+\int_0^T\|\tilde z(t)\|_{\tilde H^{\frac 32}}^4dt
\right].
\]
With these bounds on the $L^\infty(0,T;\tilde H^1)$ and
$L^2(0,T;\tilde H^2)$-norms, we can get estimates  for the time
derivative, i.e.
$\frac {d\tilde v^\lambda}{dt}\in L^2(0,T;\tilde H^0)$, in order to
conclude the proof as before.

Since $\tilde v^\lambda\in C([0,T];\tilde H^1)\cap L^2(0,T;\tilde H^2)
\subset L^4(0,T;\tilde H^{\frac 32})$ and by Proposition \ref{p:ou}
 $\tilde z \in C([0,T];\tilde H^{\frac 32})$, we get that 
$\tilde u^\lambda\in C([0,T];\tilde H^1)\cap L^4(0,T;\tilde H^{\frac 32})$.

The continuous dependence on the initial data is obtained as in the
previous proposition. Let as consider two solutions with different
initial data and the equation \eqref{eq:diff} for the difference.
Then, with usual procedure and using Lemma \ref{lemma su B tilde} 
\[\begin{split}
\frac 12 \frac {d}{dt}& \|U^\lambda(t)\|_{\tilde H^1}^2
+\nu  \|U^\lambda(t)\|_{\tilde H^2}^2
\\&=
-\langle \tilde B^\lambda(U^{\lambda}(t),\tilde u^{\lambda,(1)}(t))+
  \tilde B^\lambda(\tilde u^{\lambda,(2)}(t),U^{\lambda}(t)), 
\tilde A U^\lambda(t)\rangle
\\
& \le 
\|\tilde B^\lambda(U^{\lambda}(t),\tilde u^{\lambda,(1)}(t))+
  \tilde B^\lambda(\tilde u^{\lambda,(2)}(t),U^{\lambda}(t))\|_{\tilde H^0} 
\|U^\lambda(t)\|_{\tilde H^2}
\\ 
&\le C (1+|\lambda|) \left(\|U^{\lambda}(t)\|_{\tilde H^{\frac 12}}
\|\tilde u^{\lambda,(1)}(t)\|_{\tilde H^{\frac 32}}+
\|\tilde u^{\lambda,(2)}(t)\|_{\tilde H^{\frac 12}}
\|U^{\lambda}(t)\|_{\tilde H^{\frac 32}}\right)
\|U^\lambda(t)\|_{\tilde H^2}
\\
&\le
C (1+|\lambda|) \left(\|\tilde u^{\lambda,(1)}(t)\|_{\tilde H^{\frac 32}}+
\|\tilde u^{\lambda,(2)}(t)\|_{\tilde H^{\frac 32}}\right)
\|U^{\lambda}(t)\|_{\tilde H^{\frac 32}}
\|U^\lambda(t)\|_{\tilde H^2}
\\
& \le
C (1+|\lambda|) \left(\|\tilde u^{\lambda,(1)}(t)\|_{\tilde H^{\frac 32}}+
\|\tilde u^{\lambda,(2)}(t)\|_{\tilde H^{\frac 32}}\right)
\|U^{\lambda}(t)\|_{\tilde H^{1}}^{\frac 12}
\|U^\lambda(t)\|_{\tilde H^2}^{\frac 32}
\\
& \le 
\frac \nu 2 \|U^\lambda(t)\|_{\tilde H^2}^2
+ C(1+\lambda^2)\left(\|\tilde u^{\lambda,(1)}(t)\|^4_{\tilde H^{\frac 32}}+
\|\tilde u^{\lambda,(2)}(t)\|^4_{\tilde H^{\frac 32}}\right)
\|U^{\lambda}(t)\|_{\tilde H^{1}}^2.
\end{split}\]
Using Gronwall lemma, with usual procedure  we get
\[
\|\tilde u^{\lambda,(1)}-\tilde u^{\lambda,(2)}\|_{C([0,T];\tilde H^1)}\le C 
\|\tilde x^{(1)}-\tilde x^{(2)}\|_{\tilde H^1}
\]
for a suitable constant depending on
$T,\lambda,\nu$ and $\|\tilde u^{\lambda,(i)}\|_{L^4(0,T;\tilde H^{\frac 32})}$, 
$i=1,2$. This gives continuous dependence on the initial data
 and thus Feller property in $\tilde H^1$.
\hfill$\Box$

\medskip
The previous results are classical and we could have skipped the details
of the proof,
quoting previous results for stochastic 2D Navier-Stokes equations.
We have given all the details in order to show that 
the technique used in previous papers 
for the  stochastic 2D Navier-Stokes equations
is successful also for our system \eqref{NSE_lambda}.
With this remark in mind, we can now state another regularity result, 
based on the technique of \cite{Ferrario99}, thanks to the
properties \eqref{perdimo1}-\eqref{perdimo2}
of the bilinear operator $\tilde B^\lambda$.
\begin{proposition} \label{proc}
Let $\alpha \in \mathbb N, \alpha \ge 2$.\\
Assume $\ra \subseteq H^{\alpha+\epsi}$ for some $\epsi >0$.
Then for arbitrary $\tilde{x} \in \tilde H^\alpha$,
there exists  a unique process $\tilde u^\lambda$ solution to system 
\eqref{NSE_lambda}  over $[0,T]$ with the initial condition
$\tilde u^\lambda(0)=\tilde{x}$, such that
\[
 \tilde u^\lambda \in C([0,T];\tilde H^\alpha) 
\] 
$\po-$a.s.; 
it is a Markov  and Feller process 
in $\tilde H^\alpha$.
\end{proposition}

We recall that $H^\alpha\subset [C(\mathcal D)]^2$ for $\alpha>1$.
Therefore, by means of the latter Proposition we can 
define (for each fixed time and space) 
the quantity $\tilde u^\lambda(t,l\hat \xi)-\tilde u^\lambda(t,0)$ appearing in the 
structure functions  (for $\lambda\neq 0$ and for
$\lambda=0$ as well).

\subsection{Existence of invariant measures}

The existence of an invariant measure is obtained by means of
Krylov-Bogoliubov method (see e.g. \cite{DPZ2}). 
This requires Feller property in $\tilde H^0$ and a tightness result.
In particular,  following Chow (see Theorem 2.2 in \cite{Chow}) 
it is enough to prove the tightness for the time averages in the
following form
\begin{equation}\label{dis-chow}
\lim_{R \to \infty}\sup_{T>T_0}
\frac 1T \int_0^T\po\{\|\tilde u^\lambda(t)\|_{\tilde H^1}>R\}dt=0.
\end{equation}
for some $T_0>0$. This provides that there exists an invariant measure
with support in $\tilde H^1$.

To prove it,  we first 
show  that the  solution $\tilde u^\lambda$ on the time interval
$[0,T]$ fulfills 
\begin{equation}\label{stima-media}
\mathbb E \|\tilde u^\lambda(t)\|_{\tilde H^0}^2+2\nu \mathbb E \int_0^t \|\tilde u^\lambda(s)\|_{\tilde H^1}^2 ds 
= \|\tilde{x}\|^2_{\tilde H^0}+  2
\|G\|_{HS}^2 \ t
\end{equation}
for any $t \in (0,T]$. 
This can be done as in \cite{VF} (Theorem 1.2, Ch X), by means of
It\^o formula for 
$d\|\tilde u^\lambda(t)\|^2_{\tilde H^0}$; indeed,  assuming  {\bf H0} and using 
$\langle \tilde{B}^{\lambda}(\tilde u^\lambda,\tilde u^\lambda),
\tilde u^\lambda\rangle=0$ for $\tilde u^\lambda \in \tilde H^1$, we get
\[
\|\tilde u^\lambda(t)\|_{\tilde H^0}^2
+2\nu  \int_0^t \|\tilde u^\lambda(s)\|_{\tilde H^1}^2 ds 
= \|\tilde{x}\|_{\tilde H^0}^2+
2\int_0^t \langle \tilde u^\lambda(s), \tilde G d\tilde \beta(s)\rangle  +
2  \|G\|_{HS}^2 \ t
\]
and we can proceed as in \cite{VF} to obtain \eqref{stima-media}.

Now, taking the solution of \eqref{NSE_lambda} on the time 
interval $[0,T]$ with $\tilde x=0$, relationship \eqref{stima-media} becomes
\[
\mathbb E \|\tilde u^\lambda(T)\|_{\tilde H^0}^2+2\nu \mathbb E \int_0^T \|\tilde u^\lambda(s)\|_{\tilde H^1}^2 ds 
= 2 
\|G\|_{HS}^2 \ T
\]
giving
\[
\frac 1T   \int_0^T \mathbb E\|\tilde u^\lambda(t)\|_{\tilde H^1}^2 dt 
\le  
\frac 1\nu \|G\|_{HS}^2 .
\]
Then, by Chebyshev inequality we have
\[
\frac 1T \int_0^T \po\{\|\tilde u^\lambda(t)\|_{\tilde H^1}>R\}dt
\le 
\frac 1T \int_0^T\frac 1{R^2} \mathbb E \|\tilde
u^\lambda(t)\|_{\tilde H^1}^2 dt \le \frac 1{R^2} \frac{\|G\|_{HS}^2}{\nu}
\]
proving \eqref{dis-chow}.

Therefore we have proven the following result.
\begin{theorem} \label{teo:es-mis}
Assume {\bf H0}.
\\
Then, there exists at least one 
invariant measure $\mu^\lambda$ for system \eqref{NSE_lambda}.
Moreover $\mu^\lambda(\tilde H^1)=1$.
\end{theorem}

%%%%%%%%%%%%%%%%%%%%%%%%%%%%%%%%%%%
\subsection{Uniqueness of invariant measures}
In this section we prove that there exists a unique invariant 
measure $\mu^\lambda$
and its support is in $[C(\mathcal D)]^2$, which is important to
define the structure functions as we explained  after Proposition 
\ref{proc}. We actually prove that $\mu^\lambda(\tilde H^2)=1$ 
under suitable assumption on the 
covariance of the noise.

Uniqueness of the invariant measure can  be proven 
by different methods (see \cite{KSb}). 
Here we follow \cite{DPZ2}: we fix $\alpha\ge 2$ (as in Proposition
\ref{proc}) and  show that the Markov semigroup
$\{P^\lambda_t\}_{t\ge 0}$ is irreducible and strong Feller in $\tilde H^\alpha$. 
By means of Khasminski and Doob theorems we get
uniqueness of the invariant measure, which is strongly mixing and
equivalent to all the transition functions.

Let us recall the definitions.
Irreducibility in $\tilde H^\alpha$ means that $P^\lambda(t, \tilde{x},\Gamma)>0$ for
any $t>0$, $\tilde{x}\in \tilde H^\alpha$ and open non-empty subset $\Gamma$
of $\tilde H^\alpha$. The Markov semigroup is strongly Feller in $\tilde H^\alpha$
if $P^\lambda_t:B_b(\tilde H^\alpha) \rightarrow C_b(\tilde H^\alpha)$ for any
$t>0 $.

We shall prove in the next two subsections the following result, for each
fixed $\lambda \in \mathbb R$ and  $\alpha \in \mathbb N, \alpha \ge 2$.
\begin{theorem} \label{globale}
Assume that the operator $G$ is injective and there exists $\epsi > 0$ such that
$$
H^{\alpha+1} \subseteq \ra \subseteq H^{\alpha+\epsi} .
$$
Then the Markov semigroup $P_t^\lambda$ is irreducible and strongly Feller in
$\tilde H^\alpha$. Therefore there exists a unique invariant measure
$\mu^\lambda$ of the equation \eqref{NSE_lambda}; it is  
supported on $\tilde H^\alpha$, it is  equivalent 
to each transition probability $P^\lambda(t,\tilde x,\cdot)$ and 
\begin{equation} \label{convmis}
\lim_{t \rightarrow +\infty}P^\lambda(t,\tilde x,\Gamma)= \mu^\lambda(\Gamma)
\qquad \forall \Gamma \in {\mathcal B}(\tilde H^\alpha)
\end{equation}
for arbitrary $\tilde x \in \tilde H^\alpha$.
\end{theorem}
We will obtain this result by merging Propositions \ref{pro-irr} and
\ref{propSF}.
Notice that if $H^{\alpha+1} \subseteq \ra \subseteq H^{\alpha+\epsi}$ 
then $\ra$ is dense in $H^{\alpha}$, since $H^{\alpha+1}$ is densely 
embedded into $H^{\alpha}$. 

For example, our assumption is fulfilled if we choose $G=A^{-p}$ with
$p>1$; in this case  $\ra=H^{2p}$ and therefore
there exists an integer $\alpha\ge 2$ for which the assumptions
 of Theorem \ref{globale} are fulfilled. This gives a full noise with 
suitable space regularity. Otherwise, one 
 could prove uniqueness of the invariant measure with a degenerate
noise as in \cite{HM06,HM11}.

\bigskip

We divide the proof into two parts; first, we consider irreducibility
and then the strong Feller property.

%%%%%%%%%%%%%%%%%%%%%%%%%%%%%%%
\subsubsection{Irreducibility}
Irreducibility in $\tilde H^\alpha$ means that, 
starting from any $\tilde x \in \tilde H^\alpha$, there is a strictly positive
probability to be at any time $t>0$  in any non-empty open 
subset of $\tilde H^\alpha$. 
It is enough to check it for any open ball
$B(\tilde y,\rho)=\{\tilde x \in \tilde H^\alpha: 
\|\tilde x-\tilde y\|_{\tilde H^\alpha}<\rho\}$.
Therefore we have  to check that
\[
P^\lambda\left(t,\tilde x , B(\tilde y,\rho)\right)\equiv 
\po\{\|\tilde u^\lambda(t;\tilde u^\lambda(0)=\tilde x)-\tilde y\|_{\tilde H^\alpha}<\rho\}>0
\]
for any $\tilde x,\tilde y \in \tilde H^\alpha$ and $t, \rho>0$.
We proceed as in \cite{FM95, Ferrario97}.

We define $\tilde u_*:[0,t]\to \tilde H^\alpha$ linking $\tilde x$ to 
$\tilde y$ as
\[
\tilde u_*(s)=\begin{cases}
e^{-s\tilde A}\tilde x & s \in [0,\frac 14 t]\\
e^{-(t-s)\tilde A}\tilde y & s \in [\frac 34 t, t]\\
\tilde u_*(\frac t4)+\dfrac{s-\frac t4}{\frac t2}\left(\tilde u_*(\frac 34 t ) -
\tilde u_*(\frac 14 t)\right)
& s \in ]\frac 14 t,\frac 34 t[
\end{cases}
\]
This belongs to $C([0,t]; \tilde H^\alpha)$. Then we consider the
(unique) solution $\tilde v_*^\lambda$ of 
\[
\frac {d\tilde v_*^\lambda}{dt}(t)+\nu \tilde A \tilde v_*^\lambda (t)=
-\tilde B^\lambda(\tilde u_*(t), \tilde u_*(t));
\qquad \tilde v_*^\lambda(0)=\tilde{x}
\]
By Lemma \ref{lemma su B tilde} the r.h.s. belongs to 
$C([0,t];\tilde H^{\alpha-1})$. Then 
$\tilde v_*^\lambda\in C([0,t];\tilde H^\alpha)\cap L^2(0,t;\tilde H^{\alpha+1})$.
Finally we set 
$\tilde z^\lambda_*=\tilde u_*-\tilde v^\lambda_*\in C([0,t];\tilde
H^\alpha)$
and thus the equation fulfilled by $\tilde v_*^\lambda$ can be written
also as
\[
\frac {d\tilde v_*^\lambda}{dt}+\nu \tilde A \tilde v_*^\lambda 
+\tilde B^\lambda(\tilde v^\lambda_*+\tilde z^\lambda_*,\tilde v^\lambda_*+\tilde z^\lambda_*)=0;
\qquad \tilde v_*^\lambda(0)=\tilde{x}
\]

Now we prove a continuous dependence of $\tilde v^\lambda$ on $\tilde
z$ in equation \eqref{eq-v}. This is a deterministic result and is
proven for any integer $\alpha\ge 2$.
\begin{lemma}
We are given $\tilde x \in \tilde H^\alpha$ and 
$\tilde z_1,\tilde z_2 \in C([0,T];\tilde H^\alpha)$.
Let $\tilde v^\lambda_i\in C([0,T];\tilde H^\alpha)$ be the solution of 
\[
\frac {d\tilde v_i^\lambda}{dt}(t)+\nu \tilde A \tilde v_i^\lambda (t)
+\tilde B^\lambda(\tilde v_i^\lambda(t)+\tilde z_i(t),\tilde
v_i^\lambda(t)+\tilde z_i(t))=0, \qquad \tilde v_i^\lambda(0)=\tilde x
\]
for $i=1,2$.
Then, there exists a constant $C$ (depending on $T, \nu, \lambda,
\|\tilde v_1^\lambda+\tilde z_1\|^2_{\tilde H^{\alpha}}$ and
$\|\tilde v_2^\lambda+\tilde z_2\|^2_{\tilde H^{\alpha}}$) such that
\[
\|\tilde v_1^\lambda-\tilde v_2^\lambda\|_{C([0,T];\tilde H^\alpha)}
\le C \|\tilde z_1-\tilde z_2\|_{C([0,T];\tilde H^\alpha)}.
\]
\end{lemma}
\proof
Set $\tilde V^\lambda=\tilde v_1^\lambda-\tilde v_2^\lambda$
and $\tilde Z=\tilde z_1-\tilde z_2$. Then $\tilde V^\lambda$
satisfies
\begin{multline}
\frac {d\tilde V^\lambda}{dt}(t)+\nu \tilde A \tilde V^\lambda(t)
+\tilde B^\lambda(\tilde V^\lambda(t)+\tilde Z(t),\tilde
v_2^\lambda(t)+\tilde z_2(t))\\
+
\tilde B^\lambda(\tilde v_1^\lambda(t)+\tilde z_1(t),\tilde
V^\lambda(t)+\tilde Z(t))=0
\end{multline}
with $\tilde V^\lambda(0)=0$. We multiply this equation by $\tilde
A^\alpha \tilde V^\lambda(t)$ and integrate on the domain:
\begin{multline}
\frac 12 \frac{d}{dt}\|\tilde V^\lambda(t)\|^2_{\tilde H^\alpha}+
\nu \|\tilde V^\lambda(t)\|^2_{\tilde H^{\alpha+1}}
=\\
\langle \tilde A^{\frac {\alpha-1}2}[\tilde B^\lambda(\tilde V^\lambda(t)+\tilde Z(t),\tilde
v_2^\lambda(t)+\tilde z_2(t))+\tilde B^\lambda(\tilde v_1^\lambda(t)+\tilde z_1(t),\tilde
V^\lambda(t)+\tilde Z(t))],
\tilde A^{\frac {\alpha+1}2}\tilde V^\lambda(t)\rangle.
\end{multline}
We estimate the bilinear terms by means of Lemma \ref{lemma su B
  tilde} iii)
\[\begin{split}
\|\tilde B^\lambda(\tilde V^\lambda&+\tilde Z,\tilde
v_2^\lambda+\tilde z_2)+\tilde B^\lambda(\tilde v_1^\lambda+\tilde z_1,\tilde
V^\lambda+\tilde Z)\|_{\tilde H^{\alpha-1}}
\\&\le C
\left[\|\tilde v_1^\lambda+\tilde z_1\|_{\tilde H^{\alpha}}
+\|\tilde v_2^\lambda+\tilde z_2\|_{\tilde H^{\alpha}}\right] 
\|V^\lambda+\tilde Z\|_{\tilde H^{\alpha}}
\end{split}
\]
so that the r.h.s. is bounded by
\[
C
\left[\|\tilde v_1^\lambda+\tilde z_1\|_{\tilde H^{\alpha}}
+\|\tilde v_2^\lambda+\tilde z_2\|_{\tilde H^{\alpha}}\right] 
\|V^\lambda+\tilde Z\|_{\tilde H^{\alpha}} \|V^\lambda\|_{\tilde H^{\alpha+1}}.
\]
By Young inequality this is bounded by
\[
\frac \nu 2 \|V^\lambda\|_{\tilde H^{\alpha+1}}^2+C_\nu 
\left[\|\tilde v_1^\lambda+\tilde z_1\|^2_{\tilde H^{\alpha}}
+\|\tilde v_2^\lambda+\tilde z_2\|^2_{\tilde H^{\alpha}}\right] 
\|V^\lambda+\tilde Z\|^2_{\tilde H^{\alpha}}.
\]
Therefore, setting $\phi=\|\tilde v_1^\lambda+\tilde z_1\|^2_{\tilde H^{\alpha}}
+\|\tilde v_2^\lambda+\tilde z_2\|^2_{\tilde H^{\alpha}}\in C([0,T])$, we
obtain 
\[
 \frac{d}{dt}\|\tilde V^\lambda(t)\|^2_{\tilde H^\alpha}
\le C \phi(t) \|V^\lambda(t)\|^2_{\tilde H^{\alpha}}
 +C \phi(t)\|\tilde Z(t)\|^2_{\tilde H^{\alpha}}.
\]
Since $\tilde V^\lambda(0)=0$, 
Gronwall lemma gives the required result. \hfill $\Box$

\medskip
Since $\tilde u^\lambda-\tilde u_*=
\tilde v^\lambda-\tilde v^\lambda_*+\tilde z-\tilde z^\lambda_*$, from
the triangle inequality and the latter lemma it
follows that there exists a constant $\tilde C$ such that
\[
\|\tilde u^\lambda-\tilde u_*\|_{C([0,t];\tilde H^\alpha)}
\le \tilde C \|\tilde z-\tilde z^\lambda_*\|_{C([0,t];\tilde H^\alpha)}.
\]

Now we come back to estimate
$
\po\{\|\tilde u^\lambda(t;\tilde u^\lambda(0)=\tilde x)-\tilde
y\|_{\tilde H^\alpha}<\rho\}
$. 
Since the inital data is always $\tilde x$, in the sequel 
we drop it for simplicity.
We have
\begin{equation}\label{le-po}\begin{split}
\{\|\tilde u^\lambda(t)-\tilde
y\|_{\tilde H^\alpha}<\rho\}
&=
\{\|\tilde u^\lambda(t)-\tilde
u_*(t)\|_{\tilde H^\alpha}<\rho\}
\\&
\supseteq 
\{\|\tilde u^\lambda-\tilde
u_*\|_{C([0,t];\tilde H^\alpha)}<\rho\}
\\&
\supseteq 
\{\|\tilde z-\tilde
z^\lambda_*\|_{C([0,t];\tilde H^\alpha)}<\frac \rho{\tilde C}\}
\end{split}\end{equation}
where $\tilde u^\lambda, \tilde z$ are processes and 
$\tilde u_*, \tilde z^\lambda_*$ are deterministic functions.
We want to show that the latter term is strictly positive.
Properties of the Ornstein-Uhlenbeck process 
$\tilde z$ have been given in Proposition \ref{p:ou}. 
If   we assume in addition 
that $\ra$ is dense in $\tilde H^\alpha$, then 
by Proposition 2.7. in \cite{Masl} we have that the closure of the support 
law of is $C_0([0,t];\tilde H^\alpha)$ (when the subscript $0$ denotes
that the initial value vanishes). Therefore the law of the process 
 $\tilde z$ is a full
measure in $C_0([0,t];\tilde H^\alpha)$, i.e.
\[
\po\{\|\tilde z-\tilde \zeta\|_{C([0,t];\tilde H^\alpha)}<r\}>0
\]
for any $\tilde \zeta \in C_0([0,t];\tilde H^\alpha)$ and $r>0$.
Keeping in mind \eqref{le-po} we obtain
\[
\po\{\|\tilde u^\lambda(t)-\tilde
y\|_{\tilde H^\alpha}<\rho\}>0.
\]

Summing up,
for any  $\lambda \in \mathbb R$ and  
$\alpha \in \mathbb N, \alpha \ge 2$  we have proved
\begin{proposition}\label{pro-irr}
Assume  $\ra \subseteq H^{\alpha+\epsi}$ for some $\epsi>0$ 
and $\ra$ dense in $H^{\alpha}$.
Given any  $\tilde x, \tilde y \in \tilde H^\alpha$, $t>0$ and
$\rho>0$  we have 
\[
P^\lambda(t,\tilde x,B(\tilde y,\rho)) > 0
\]
i.e. the Markov semigroup
$P_t^\lambda$ is irreducible in $\tilde H^\alpha$.
\end{proposition}

%%%%%%%%%%%%%%%%%%%%%%%%%%%%%%%
\subsubsection{Strong Feller}
The second property of the Markov semigroup we have to check
is the strongly Feller property in $\tilde H^\alpha$, i.e.
$P^\lambda_t:B_b(\tilde H^\alpha) \rightarrow C_b(\tilde H^\alpha)$ for any
$t>0$.\\
We already proved that  $\{P^\lambda_t\}_{t\ge 0}$ 
is Feller in $\tilde H^\alpha$. 
 By the mean value Theorem, we would get that it is
Lipschitz Feller if we were able to estimate the derivative of
$P_t^\lambda \psi$. This is not true, but  as in
\cite{Ferrario99} we can prove  it for a
modified version of \eqref{NSE_lambda}
\begin{multline} \label{ns-tronc} 
        d\tilde u^{\lambda,(R)}(t) +
        \tilde A\tilde u^{\lambda,(R)}(t)dt\\
       +\Theta_{R}(\|\tilde u^{\lambda,(R)}(t)\|^{2}_{\tilde H^\alpha})
        \tilde B^\lambda(\tilde u^{\lambda,(R)}(t),\tilde u^{\lambda,(R)}(t))\ dt
    =\tilde G  d\tilde \beta(t)
\end{multline}
where cut-off function
$\Theta_{R}$ is a $C^\infty$
function equal to 1 in $[-R,R]$ and 0 outside $[-R-1,R+1]$.

By means of Bismut-Elworthy-Li's formula, we prove that
the Markov semigroup associated to \eqref{ns-tronc} is Lipschitz
Feller in $\tilde H^\alpha$.
\begin{proposition}\label{pro:sfR}
Assume that for some $\alpha \in \mathbb N$ with $\alpha\ge 2$ the
operator 
$G$ is injective with
$H^{\alpha+1} \subseteq \ra \subseteq H^{\alpha+\epsi}$
for some $\epsi > 0$.\\
Then,  for every $\lambda \in \mathbb R$ and $t,R>0$
there exists a constant $L=L(\lambda,R,t)$ such that
\[
\left|P_{t}^{\lambda, (R)} \psi (\tilde x) 
  - P_{t}^{\lambda, (R)} \psi(\tilde y) \right| \; \le
 \; L \;
\|\tilde x-\tilde y\|_{\tilde H^\alpha}
\]
for all $\tilde x,\tilde y \in \tilde H^\alpha,\psi \in C_b(\tilde H^\alpha)$
with $\|\psi\|_b \le 1$.\\
Moreover, $P_{t}^{\lambda, (R)} \psi$ is Lipschitz continuous for arbitrary
$\psi \in B_b(\tilde H^\alpha)$.
\end{proposition}

The proof is the same as in \cite{FM95,Ferrario99}, based on the properties
of the bilinear operator $\tilde B^\lambda$ given in Lemma \ref{lemma
  su B tilde}. 

Moreover, from the estimates of Section 3.1 and Proposition 3.3, we have that 
 $u^\lambda, u^{\lambda, (R)}\in C([0, T]; \tilde H^\alpha)$
$\po$-a.s. and 
\[
\sup_{\|\tilde{x}\|_{\tilde H^\alpha}\le M}\sup_{0\le t\le T}
\|u^{\lambda, (R)}(t, \tilde{x})\|_{\tilde H^\alpha}<\infty,
\]
Thus the processes $u^\lambda(\cdot, \tilde{x})$ and  
$u^{\lambda, (R)}(\cdot, \tilde{x})$ coincide in the ball 
$B_{R}:=\{v:\ \|v\|^{2}_{\tilde H^\alpha}\le R\}$.
Therefore one  proves that 
\begin{equation}\label{eq:limR}
\lim_{R \to \infty}\|P^{\lambda,(R)}(t,\tilde x,\cdot)-
P^{\lambda}(t,\tilde x,\cdot)\|_{var}=0
\end{equation}
uniformly with respect to $\tilde x$ in bounded sets of $\tilde
H^\alpha$,
where $\|\cdot\|_{var}$ denotes the total variation norm of a measure.

Now, passing to the limit as  $R \rightarrow \infty$,
we  obtain the strong
Feller property for the principal equation \eqref{NSE_lambda}.
Indeed, 
\begin{multline*}
\|P^{\lambda}(t,\tilde x,\cdot)-
P^{\lambda}(t,\tilde y,\cdot)\|_{var}\le
\|P^{\lambda}(t,\tilde x,\cdot)-
P^{\lambda,(R)}(t,\tilde x,\cdot)\|_{var}
\\+
\|P^{\lambda,(R)}(t,\tilde x,\cdot)-
P^{\lambda,(R)}(t,\tilde y,\cdot)\|_{var}+
\|P^{\lambda,(R)}(t,\tilde y,\cdot)-
P^{\lambda}(t,\tilde y,\cdot)\|_{var}
\end{multline*}

We fix any $\epsi>0$.
From \eqref{eq:limR}, there exists $R_\epsi>0$ such that
\[
\|P^{\lambda}(t,\tilde x,\cdot)-
P^{\lambda,(R_\epsi)}(t,\tilde x,\cdot)\|_{var}<\epsi ,
\]
\[
\|P^{\lambda,(R_\epsi)}(t,\tilde y,\cdot)-
P^{\lambda}(t,\tilde y,\cdot)\|_{var}<\epsi .
\]
On the other hand,
from Proposition \ref{pro:sfR}
there exists  $\delta_\epsi>0$ such that 
\[
\|P^{\lambda,(R_\epsi)}(t,\tilde x,\cdot)-
P^{\lambda,(R_\epsi)}(t,\tilde y,\cdot)\|_{var}
< \epsi
\]
for all $\tilde x, \tilde y$ with 
$\|\tilde x-\tilde y\|_{\tilde H^\alpha}<\delta_\epsi$. 

Thus, given any $\epsi>0$ there exists  $\delta_\epsi>0$ such that 
\[
\|P^{\lambda}(t,\tilde x,\cdot)-
P^{\lambda}(t,\tilde y,\cdot)\|_{var}
< 3\epsi
\]
for all $\tilde x, \tilde y$ with 
$\|\tilde x-\tilde y\|_{\tilde H^\alpha}<\delta_\epsi$. 
Therefore we have proven the strong Feller property.
\begin{proposition}\label{propSF}
Assume that for some $\alpha \in \mathbb N$ with $\alpha\ge 2$ the
operator 
$G$ is injective and 
$H^{\alpha+1} \subseteq \ra \subseteq H^{\alpha+\epsi}$
for some $\epsi > 0$.\\
Then, $P^\lambda_t$ is strong Feller in $\tilde H^\alpha$.
\end{proposition}

%%%%%%%%%%%%%%%%%%%%%%%%%%%%%%%
\section{Continuous dependence of the invariant measure wrt 
to the parameter $\lambda$}\label{conv}
In this section we show the continuous dependence of the invariant
measure $\mu^\lambda$ on the parameter $\lambda$. In particular we are
interested in the case of $\lambda \to 0$.

We have a first result.
Let us  fix the family of the unique invariant measures $\mu^\lambda$ (given in
Section 3), and consider the limit when $\lambda\to 0$.

\begin{proposition} \label{tightness}
The family of  
invariant measures $\{\mu^{\lambda}\}_{\lambda\neq 0}$ is tight in
$\tilde H^{1-\epsilon}$ for
any $\epsilon>0$.  Therefore there exists a measure $m$ on
$H^{1-\epsilon}$ and a sequence $\{\mu^{\lambda_n}\}_{n \in \mathbb  N}$ 
(with $\lambda_n \to 0$ as $n \to \infty$) weakly converging 
 to $m$ in $ \tilde H^{1-\epsilon}$, i.e.
\[
\lim_{n \to \infty}\int f d\mu^{\lambda_n} =\int f dm\qquad
\forall f \in C_b(\tilde H^{1-\epsilon}).
\]
Finally, the supports of $m$ and $\mu^{\lambda_n}$ are contained in $\tilde H^1$.
\end{proposition} 
\proof 
 For each $\lambda \in \mathbb  R$ 
let us denote by $\tilde u^\lambda_{st} $ the 
stationary process solutions of system 
\eqref{NSE_lambda} whose law at each fixed time is $\mu^\lambda$, 
the invariant measure defined in Section 3.
From \eqref{stima-media} we get that
\begin{equation}\label{ito-invariant-1}
\nu {\mathbb E} \int_{0}^{t} \|\tilde u_{st}^\lambda(s)\|_{\tilde H^1}^2 ds 
= \|G\|_{HS}^2 \ t
\end{equation}
i.e.
\begin{equation}\label{ito-invariant-12}
 \int \| x\|^2_{\tilde H^1} d\mu^\lambda(x) =  \frac{\|G\|_{HS}^2} {\nu}
\end{equation}
uniformly in $\lambda$.
Then, using that $\tilde H^1$ is compactly
embedded in $\tilde H^{1-\epsilon}$ we get tightness by means of  the Chebyshev 
inequality. Hence, by Prohorov theorem 
there exists a subsequence $\mu^{\lambda_n}$ that converges to 
a probability measure $m$ in $\tilde H^{1-\epsilon}$. 

As far as the supports are concerned, from Theorem \ref{teo:es-mis}
we know that $\mu^\lambda(\tilde
H^1)=1$ for each $\lambda$. According to Lemma 3.1 in Ch II of
\cite{VF}, from   \eqref{ito-invariant-1} we get that also the limit
measure $m$ has support contained in $\tilde H^1$.
\hfill $\Box$

\medskip
Now we have to show that the limit measure $m$ is the unique invariant
 measure associated to system \eqref{Joint}.
For this purpose, we have to work with the stationary process solving
\eqref{NSE_lambda}. 
In the next subsection we shall prove the following result
\begin{proposition}\label{proc-tight}
The family $\{\tilde u^\lambda_{st}\}_{\lambda \in \mathbb R}$ of stationary processes solving 
\eqref{NSE_lambda} is tight in 
$L^2_{loc}(0,\infty;\tilde H^{0})\cap C([0,\infty); \tilde H ^{-1})$.

Thus there exists a new probability  basis $(\bar \Omega, \bar  F, \bar P)$,
a sequence $\{\bar u_{st}^{\lambda_n}\}$ of stationary processes and a
limit   process $\bar u_{st}^0$ 
defined on it with values in $L^2_{loc}(0,\infty;\tilde H^{0})\cap
C([0,\infty); \tilde H ^{-1})$ and solving \eqref{NSE_lambda} with
  parameter $\lambda_n$ and 0 respectively,
such that $\tilde u^{\lambda_n}_{st}$ and $\bar u_{st}^{\lambda_n}$ have the same law and 
\[
\lim_{n \to \infty}\bar u_{st}^{\lambda_n}=\bar u_{st}^0 \text{ in }
L^2_{loc}(0,\infty;\tilde H^{0})\cap C([0,\infty); \tilde H ^{-1}) \qquad \bar P-a.s.
\]
Finally, the process $\bar u_{st}^0$ is a stationary process in $\tilde H^0$.
\end{proposition}
In particular, for any $t\ge 0$ we have
\[
\lim_{n \to \infty}\bar u_{st}^{\lambda_n}(t)=\bar u_{st}^0(t) \text{ in }
\tilde H^{-1} \qquad \bar P-a.s.
\]

Since the law of $\bar u_{st}^{\lambda_n}(t)$ is $\mu^{\lambda_n}$ and the
law of $\bar u_{st}^0(t)$ is $\mu^0$, we have that
\[
\lim_{n \to \infty}\mu^{\lambda_n}= \mu^0 
 \quad {\rm weakly \ in }\quad \tilde H^{-1} .
\]
Since there exists a unique invariant measure for the system \eqref{Joint}
with $\lambda=0$, we get that any sequence extracted from 
$\{\mu^\lambda\}_{\lambda \in \mathbb R}$ weakly converges to $\mu^0$
in $\tilde H^{-1}$ as $\lambda\to 0$. Bearing in mind Proposition 
\ref{tightness}, this identifies $m$ with $\mu^0$ as measures on 
Borelian sets of $\tilde H^{-1}$. Since we know that 
both $m$ and $\mu^0$ are supported on $\tilde H^1$ indeed, we get that 
 $m=\mu^0$ as measures on 
Borelian subsets of $\tilde H^{1}$.
\begin{theorem}
Let $\epsi>0$ be given.
For any sequence  $\{\mu^{\lambda_n}\}_{n \in \mathbb  N}$ 
(with $\lambda_n \to 0$ as $n \to \infty$) we have
\[
\lim_{n \to \infty}\mu^{\lambda_n} =  \mu^0 
 \quad {\rm weakly \ in }\quad \tilde H^{1-\epsi} .
\]
\end{theorem}

%%%%%%%%%%%%%%%%%%%%%%%%%%%
\subsection{Convergence of stationary solutions}
We now prove Proposition \ref{proc-tight}. For simplicity we drop the
subindex and denote by $u^\lambda$ the stationary solution of \eqref{NSE_lambda} whose
marginal at any fixed time is the unique invariant measure $\mu^\lambda$.

The proof is based on two steps: first we show that the 
sequence of laws of $\tilde u^\lambda$, $\lambda >0$, is tight; 
then we pass to the limit in a suitable way and get that
the limit process is a weak solution of system \eqref{Joint} (in the probabilistic sense).

Let us recall some of the estimates performed in Section 3 
by means of  the It\^o formula: for $t\ge 0$
\begin{equation}\label{ito-formula}
\|\tilde u^\lambda(t)\|_{\tilde H^0}^2
+2\nu  \int_0^t \|\tilde u^\lambda(s)\|_{\tilde H^1}^2 ds 
= \|\tilde u^\lambda(0)\|_{\tilde H^0}^2+
2\int_0^t \langle \tilde u^\lambda(s), \tilde G d\tilde \beta(s)\rangle  +
2  \|G\|_{HS}^2 \ t.
\end{equation}
Now, using the Burkholder-Davis-Gundy inequality and taking the 
expected values yields a uniform estimate with respect to 
 $\lambda$, that is   
\[
 \sup_{\lambda\in \mathbb R}\left[\mathbb E \sup_{0\le t\le T}
\|\tilde u^\lambda(t)\|_{\tilde H^0}^2
+\nu \mathbb E \int_0^T \|\tilde u^\lambda(s)\|_{\tilde H^1}^2 ds \right]
\le  C.
\]
%where the constant $C$ is independent of $\lambda$.

Now, we write equation \eqref{NSE_lambda} in the integral form 
$$
 \tilde u^\lambda(t)=\tilde u^\lambda(0)
-\int_0^t [  \nu\widetilde{A}\widetilde{u}^{\lambda}(s)
 +\widetilde{B}^{\lambda}\left(  \widetilde{u}^{\lambda}(s),\widetilde
 {u}^{\lambda}(s)\right)  ]  ds + \tilde G \widetilde{\beta}(t),\qquad t>0.
$$
We have from  Lemma \ref{lemma su B tilde} and using an interpolation estimate
\[\begin{split}
\int_0^T\|\tilde B^\lambda(\tilde{u}^{\lambda}(s),\tilde
{u}^{\lambda}(s))\|^2_{\tilde H^{-1}} ds
&\le C^2(1+|\lambda|)^2\int_0^T \|\tilde{u}^{\lambda}(s)\|^4_{\tilde H^{\frac 12}}ds
\\&
\le  C^2 (1+|\lambda|)^2 \int_0^T \|\tilde{u}^{\lambda}(s)\|^2_{\tilde  H^0}
   \|\tilde{u}^{\lambda}(s)\|^2_{\tilde H^1}ds
\\&
\le C^2 (1+|\lambda|)^2 \|\tilde{u}^{\lambda}\|_{L^\infty(0,T;\tilde H^0)}^2
\|\tilde{u}^{\lambda}\|_{L^2(0,T;\tilde H^1)}^2
\end{split}\]
Therefore, by  usual estimations (see, e.g., \cite{FG95})
we get that there exist constants $C$ such that
\[
\begin{split}
& \sup_{\lambda\in \mathbb R}\mathbb E
   \|\int_0^\cdot \widetilde{A}\widetilde{u}^{\lambda}(s) ds\|^2_{W^{1,2}(0,T;\tilde H^{-1})}\le C\\
&\sup_{|\lambda| \le 1} \mathbb E\|\int_0^\cdot \widetilde{B}^{\lambda}\left(  \widetilde{u}^{\lambda}(s),\widetilde
 {u}^{\lambda}(s)\right) ds\|_{W^{1,2}(0,T;\tilde H^{-1})}\le C %\quad \text{ by}
 \\
&\mathbb E\| \tilde G\widetilde{\beta}(t)\|^2_{W^{\alpha,2}(0,T;\tilde H^0)}\le C(\alpha) 
     %\quad \text{ by }\eqref{noise2}
\end{split}
\]
for all $\alpha\in (0,\frac 12)$.
Therefore, for any finite $T$
$$
 \sup_{|\lambda|\le 1} \mathbb E \|\tilde u^\lambda\|_{W^{\alpha,2}(0,T;\tilde H^{-1})}<\infty.
$$
On the other hand, we already know from  the previous estimate above that
$$
 \sup_{\lambda}  \mathbb E \|\tilde u^\lambda\|^2_{L^2(0,T;\tilde H^1)}<\infty.
$$
Using that the space $L^2(0,T;\tilde H^1)\cap W^{\alpha,2}(0,T;\tilde H^{-1})$ 
is compactly embedded in $L^2(0,T; \tilde H^{0})$ and in $C([0,T]; \tilde H^{-1})$
 (see \cite{VF} Ch IV, Theorem 4.1),
it follows that the sequence of laws of processes $\left\{u^\lambda\right\}_\lambda$
is tight in $L^2(0,T; \tilde H^{0})\cap C([0,T]; \tilde H^{-1})$.

With the usual procedure (see, e.g., \cite{FG95}, we get 
that the sequence of laws of processes $\left\{u^\lambda\right\}_\lambda$
is tight in $L^2_{loc}(0,\infty; \tilde H^{0})\cap C([0,\infty); \tilde H^{-1})$.

%Let us endow $L^2_{loc}(0,\infty;H^{-s})$ by the distance
%$$
 %d_2(\xi,\zeta)=\sum_{n=1}^\infty 2^{-n} \min (\|\xi-\zeta\|_{L^2(0,n;H^{-s})},1)
%$$
%and $C([0,\infty);H^{-2-2s})$ by the distance
%$$
 %d_\infty(\xi,\zeta)=\sum_{n=1}^\infty 2^{-n} 
 %\min (\|\xi-\zeta\|_{C([0,n];H^{-2-2s})},1).
%$$
From Prokhorov and Skorohod theorems  follows that 
there exists a basis $(\bar \Omega, \bar  F, \bar P)$ (with expectation $\bar{\mathbb E}$),
and on this basis, 
$L^2_{loc}(0,\infty; \tilde H^{0})\cap C([0,\infty); \tilde H^{-1})$-valued random variables 
$\bar u^0$, $\bar u ^\lambda$, such that
${\mathcal L}({\bar u}^\lambda)={\mathcal L}({\tilde u}^\lambda)$  and
for each finite $T$
\begin{equation}\label{conv-qo}
 \lim_{\lambda_n \to 0} {\bar u}^{\lambda} = {\bar u}^0
 \qquad\text{ in }\quad L^2(0,T; \tilde H^{0})\cap C([0,T]; \tilde H^{-1})\; \bar P-a.s. 
\end{equation}
Moreover, each process $\bar u^{\lambda}$ satisfies  the same estimates
as $\tilde u^{\lambda}$ since they have the same law; hence
\[
 \sup_{\lambda\in \mathbb R}\left[\bar {\mathbb E} 
\sup_{0\leq t\le T}\|\bar u^\lambda(t)\|_{\tilde H^0}^2
+\nu \bar{\mathbb E} \int_0^T \|\bar u^\lambda(s)\|_{\tilde H^1}^2 ds \right]
\le  C.
\]
The fact that the limit process $\bar u^0$ solves system \eqref{Joint} 
follows by passing to the limit on the system \eqref{NSE_lambda}, 
see \cite{FG95} and  \cite{BF2014}.
 
In addition, $\bar u^\lambda$ and $\tilde u^\lambda$ have the same law; then 
$\bar u^\lambda$ is a stationary process. By the convergence 
$\bar P$-a.s. in $C([0,\infty); \tilde H^{-1})$ 
we get that also $\bar u^0$ is a stationary process in $\tilde H^{-1}$.

Finally,  from the estimate above,  we have that 
$$
 {\bar u}^0 \in  L^\infty(0,T; \tilde H^{0}) \qquad \bar P-a.s. 
$$
Then, for $T<\infty$ almost each path $\bar u^0 \in C([0,T];\tilde H^{-1})\cap
L^{\infty}(0,T; \tilde H^0)$; thus it is weakly continuous in $\tilde H^{0}$,
i.e. we have  for any $\phi \in \tilde H^{0}$
$$
 \lim_{t\to t_0} \int_D\bar u^0(t)\phi\ dx = 
 \int_D \bar u^0(t_0) \phi\ dx \qquad \bar P-a.s.
$$ 
and for any $t \in [0,T]$
$$
 \|\bar u^0(t)\|_{\tilde H^{0} } \le  \|\bar u^0\|_{L^\infty(0,T; \tilde H^{0} )}
 \qquad \bar P-a.s.
$$
(see \cite{Temam} p 263). \\
Hence, for every $t\ge 0$, the mapping 
$\bar \omega\mapsto \bar u^0(t,\bar\omega)$ is well
defined from $\bar \Omega$ to $\tilde H^{0}$ and it is weakly measurable. 
Since $\tilde H^{0}$
is a separable Banach space, it is strongly measurable (see
\cite{yosida} p 131). Therefore, it is meaningful to speak about the
law of $\bar u^0(t)$ in $\tilde H^{0}$. 
The stationarity of $\bar u^0$ in $\tilde H^{0}$ 
has to be understood in this sense.

%%%%%%%%%%%%%%%%%%%%%%%%%%%%%%%%
\section{Conclusions}
In this paper, we investigated the  statistics of a nonlinear model,
the stochastic Navier-Stokes system \eqref{NSE} versus its linear
counterpart given by the stochastic passive scalar equation
\eqref{Passive}. We coupled them by introducing a parameter
$\lambda\in\mathbb{R}$ and obtained a joint system
\eqref{Joint_lambda}. After rescaling the joint system, we can get  a
symmetric  system \eqref{coppia-u-v}. Moreover the system being
symmetric implies that the averages computed on each component of the
system are the same. These averages are computed with respect
to the invariant measure of the system.  

The main goal of the paper was to study the existence, uniqueness of
invariant measures for system \eqref{Joint_lambda} and their
properties with respect to the parameter $\lambda$ in particular its 
continuous dependence when $\lambda\to 0$.

We proved that the joint system \eqref{Joint_lambda} has a unique, 
ergodic invariant measure $\mu^\lambda$ for any $\lambda\in\mathbb{R}$. 
Then, when $\lambda\to 0$, we proved that  $\mu^\lambda\to \mu^0$ where 
$\mu^0$ is the unique invariant measure of joint system  
\eqref{Joint_lambda} for $\lambda=0$ which is the joint system 
\eqref{Joint}. As a consequence, the statistical properties 
obtained for the symmetric system \eqref{coppia-u-v} translate 
to the joint system \eqref{Joint}. More precisely  the statistical 
properties of \eqref{NSE} are similar to \eqref{Passive} that are 
simpler to compute.
All our results are given for a non degenerate noise but can be 
extended for a degenerate noise. 
%
%Now, if  the parameter $\lambda$ appears in the equation through the 
%noise $\tilde \beta$ instead of appearing directly in the equation 
%then a  different technique can be found in \cite{KSb} (see Theorem 4.3.1).
%
%\vspace{0.5cm}
% 
%They  proved  weak convergence of the invariant
%measures (continuous dependence wrt $\lambda$ of the invariant measures $\mu^\lambda$).
%However, Theorem 4.3.1 of \cite{KSb} requires the noise to be full, whereas our technique
%works also for a degenerate noise term. Even if we considered the full
%noise, the results of Section \ref{conv} holds  also for a degenerate
%noise.
\\[2mm]
\textbf{Acknowledgment:} Hakima Bessaih was partially supported by NSF 
grant DMS-1418838.

\end{document}